\documentclass[11pt,reqno,a4paper]{amsart}
\usepackage{amssymb,amsfonts,amsmath,amsthm}
\usepackage[dvipsnames]{xcolor}
\usepackage{graphicx}
\usepackage{lineno}
\usepackage[margin=3cm]{geometry}

\usepackage{hyperref}
\hypersetup{
    colorlinks=true,
    linkcolor=blue,
    citecolor=red,
}
\newcommand{\veps}{\varepsilon}

\newcommand{\cm}{\complement} 

\newcommand{\C}{\mathbb{C}}
\newcommand{\N}{\mathbb{N}}

\newcommand{\F}{\mathcal{F}}

\newcommand{\uld}{\overline{\operatorname{logdens}}}
\newcommand{\lld}{\underline{\operatorname{logdens}}}

\newcommand{\mes}{\operatorname{m}}

\def\r{\varrho}
\newtheorem{theorem}{Theorem}[section]

\theoremstyle{plain}

\newtheorem{corollary}[theorem]{Corollary}

{\theoremstyle{definition}

}

\newtheorem{lemma}{Lemma}[section]

\numberwithin{equation}{section}

\newtheorem{thm}{Theorem}

\title{On Petrenko's deviations and Second order differential equations}
\author{J.~Heittokangas* and M.~A.~Zemirni}

\thanks{*Corresponding author}
\begin{document}
\maketitle

\begin{abstract}
New results on the oscillation of solutions of $f''+A(z)f=0$ and on the growth of solutions of
$f''+A(z)f'+B(z)f=0$ are obtained, where $A$ and $B$ are entire functions. 
Petrenko's magnitudes of deviation of $g$ with respect to $\infty$ play a key r\^ole in the results,
where $g$ represents one of the coefficients $A$ or $B$. These quantities are defined by 
$\beta^-(\infty,g) = \liminf_{r\to\infty} \frac{\log M(r,g)}{T(r,g)}$ and
$\beta^+(\infty,g) = \limsup_{r\to\infty} \frac{\log M(r,g)}{T(r,g)}$.


\end{abstract}


\section{Introduction}

We consider the oscillation of solutions of
	\begin{equation}\label{ode}
	f'' + A(z)f=0
	\end{equation}
and the growth of solutions of
	\begin{equation}\label{ode1}
	f''+A(z)f'+B(z)f=0
	\end{equation}
where $A$ and $B$ are entire functions. It is well known that all solutions of either equation
are entire. If $g$ represents either of the coefficients $A$ or $B$, our results on the equations 
\eqref{ode} and 
\eqref{ode1} rely on the {\it magnitudes of deviation of $g$ with respect to $\infty$}
introduced by Petrenko \cite{P}. These quantities are given by
	\begin{equation}\label{deviation}
	\beta^-(\infty,g) = \liminf_{r\to\infty} \frac{\log M(r,g)}{T(r,g)}
	\quad\textnormal{and}\quad 
	\beta^+(\infty,g) = \limsup_{r\to\infty} \frac{\log M(r,g)}{T(r,g)},
	\end{equation}
where $M(r,g)=\max_{|z|=r}|g(z)|$ and $T(r,g)$ is the Nevanlinna characteristic of $g$.  

Recall that any entire function $g$ satisfies the inequalities
	\begin{equation}\label{basic-ineq}
	T(r,g)\leq \log M(r,g)\leq \frac{R+r}{R-r}T(R,g)
	\end{equation}
for all $0<r<R<\infty$ \cite[Theorem~1.6]{Hayman}. Choosing $R=2r$,
it is easy to obtain the following conclusion \cite[Theorem~1.7]{Hayman}: \emph{The functions $T(r,g)$
and $\log M(r,g)$ have the same order~$\rho$. Moreover, if $0<\rho<\infty$, then  
$T(r,g)$ and $\log M(r,g)$ are simultaneously of minimal type, mean type or maximal type.}

As for the quantities in \eqref{deviation}, if $g$ is of finite lower order $\mu$, 
then \cite[Theorem 1]{P} shows that
	\begin{equation}\label{petrenko}
	1 \le \beta^-(\infty, g) \le \mathcal{B}(\mu),
	\end{equation}
where 
	$$
	\mathcal{B}(\mu) := \left\{ 
	\begin{array}{rll}
	\displaystyle\frac{\pi \mu}{\sin(\pi \mu)}, & \ \text{if} & 0\le \mu  \le \frac{1}{2}, \\
	\pi \mu,  & \ \text{if} & \mu  \ge \frac{1}{2}.
	\end{array}
	\right.
	$$
Both of the inequalities in \eqref{petrenko} are sharp -- see \cite{C} regarding the first inequality, 
and \cite[\S 4]{P} regarding the second inequality. The construction in \cite[\S 4]{P} also
shows that $\beta^-(\infty, g)<\beta^+(\infty, g)$ may happen. The case $1< \beta^-(\infty,g) < \mathcal{B}(\mu)$ is also possible -- for example, the Airy integral $\operatorname{Ai}(z)$ has lower order $3/2$ and satisfies \cite{G0, G2}
	$$
	1 < \lim_{r\to\infty} \frac{\log M(r,\operatorname{Ai})}{T(r,\operatorname{Ai})} = \frac{3\pi}{4} < \frac{3\pi}{2}.
	$$
If $g$ is of infinite order, then $\beta^-(\infty,g)$ need not be finite. For example, if 
$g(x)=\exp\left(e^z\right)$, then $T(r,f)\sim e^r(2\pi r)^{-1/2}$ and $\log M(r,f)=e^r$ \cite[pp.~19--20]{Hayman}. 	

If $\phi(r)$ is any increasing function and convex in $\log r$ such that $\phi(r)/\log r\to \infty$,
then there is an entire function $g$ satisfying $T(r,g)\sim \phi(r)\sim \log M(r,g)$ \cite{C}. This result
allows us to construct many examples of functions $g$ for which $\beta^-(\infty, g)=\beta^+(\infty, g)$. 	

An entire function $g(z)=\sum_{n=0}^\infty a_nz^{\lambda_n}$ is said to have Fej\'er gaps if $\sum \lambda_n^{-1}<\infty$ and Fabry gaps if $\lim \lambda_n/n=\infty$. A function $g$ with Fej\'er gaps has no finite deficient values, and satisfies
	\begin{equation}\label{one}
	T(r,g)\sim \log M(r,g)
	\end{equation}
as $r\to\infty$ outside a set of finite logarithmic measure. A function $g$ with Fabry gaps
satisfies 
	\begin{equation}\label{Fabry}
	\log L(r,g)\sim \log M(r,g),\quad L(r,g)=\min_{|z|=r}|g(z)|, 	
	\end{equation}
as $r\to\infty$ outside a set of zero logarithmic density \cite{F}. Consequently, $g$ satisfies \eqref{one}
as $r\to\infty$ outside a set of zero logarithmic density. Value distribution of entire functions
$g$ satisfying \eqref{one} as $r\to\infty$ on a set of positive density is studied in \cite{HR}.

Recall that the density and the lower density of a set $F\subset [1,\infty)$ are respectively 
	$$
	\overline{\operatorname{dens}}(F)
	=\limsup_{r\to\infty}\frac{\int_{F\cap [1,r]}dt}{r-1}
	\quad\textnormal{and}\quad
	\underline{\operatorname{dens}}(F)
	=\liminf_{r\to\infty}\frac{\int_{F\cap [1,r]}dt}{r-1},
	$$
while the corresponding logarithmic densities are
	$$
	\overline{\operatorname{logdens}}(F)
	=\limsup_{r\to\infty}\frac{\int_{F\cap [1,r]}\frac{dt}{t}}{\log r}
	\quad\textnormal{and}\quad
	\underline{\operatorname{logdens}}(F)
	=\limsup_{r\to\infty}\frac{\int_{F\cap [1,r]}\frac{dt}{t}}{\log r}.
	$$
The following inequalities are known: 
	\begin{equation*}
	0\leq\underline{\operatorname{ dens}}(F)\leq\underline{\operatorname{\log dens}}(F)\leq\overline{\operatorname{\log dens}}(F)\leq\overline{\operatorname{dens}}(F)\leq 1.
	\end{equation*}

If $\beta^-(\infty,g)=\beta^+(\infty,g)<\infty$, then there exists an $\alpha\in (0,1]$ such that
	\begin{equation}\label{alpha}
	T(r,g)\sim \alpha\log M(r,g)
	\end{equation}
as $r\to\infty$ without an exceptional set. Allowing exceptional sets, this generalizes \eqref{one}. For example, $g(z)=e^z$ 
satisfies \eqref{alpha} for $\alpha=1/\pi$, while exponential polynomials in general satisfy a 
condition of the form \eqref{alpha} as $r\to\infty$ outside of a set of zero 
density \cite{LHY}.

The main results on the equations \eqref{ode} and \eqref{ode1} in terms of Petrenko's magnitudes
of deviation are stated and discussed in Sections~\ref{oscillation-sec} and \ref{growth-sec}, while the proofs of the main results are given in Sections~\ref{proof1-sec} and \ref{proof2-sec}.


\section{Oscillation theory}\label{oscillation-sec}

Bank and Laine \cite{BL,BL2} have proved that if $A$ is a transcendental entire function of order $\rho(A)< \frac{1}{2}$, then $\lambda(E)=\infty$, where $E$ is a product of two linearly independent solutions of the equation \eqref{ode} and $\lambda(g)$ denotes the exponent of convergence of zeros of $g$.
Moreover, if $\rho(A)\notin \mathbb{N}$, then
	\begin{equation}\label{bl-inequality}
	\lambda(E) \ge \rho(A).
	\end{equation}
For the case $1/2 \le \rho(A)<1$, Rossi \cite{R} improved the inequality \eqref{bl-inequality} to
	\begin{equation}\label{rs-inequality}
	\lambda(E) \ge \frac{\rho(A)}{2\rho(A)-1},
	\end{equation}
where $\lambda(E)=\infty$ if $\rho(A)=1/2$. The case $\rho(A)=1/2$ was proved independently by
Shen \cite{Shen}. Recently, Bergweiler and Eremenko have showed that \eqref{rs-inequality} is the best possible in the case $1/2 < \rho(A) < 1$ \cite{BE2}, and that \eqref{bl-inequality} is the best possible in the case $\rho(A)\ge 1$ \cite{BE}. However, under additional assumptions on the coefficient $A$, the inequality \eqref{bl-inequality} can be improved  to
	\begin{equation}\label{be-inequality}
	\lambda(E) \ge \frac{N\mu(A)}{2\mu(A)-N},
	\end{equation}
where $N$ is the number of the unbounded compenents of the set $\big\{z\in \C : |A(z)|>K|z|^p \big\}$, where $K>0$ and $p>0$, and the lower order $\mu(A)$ of $A$ satisfies $N/2 \le \mu(A)<N$ \cite{BE2}. In the same paper \cite{BE2}, Bergweiler and Eremenko showed that \eqref{be-inequality} is the best possible.

We note that, independently on the conditions imposed for the coefficient $A$, the lower bound for $\lambda(E)$ in the results stated above always depend on either $\rho(A)$ or $\mu(A)$. We proceed
to search for conditions on $A$ such that the lower bounds for $\lambda(E)$ are independent on $\rho(A)$ or $\mu(A)$. The following result by Laine and Wu is in this direction.

\begin{thm}[\cite{LW}]\label{laine-wu}
Let $A$ be a transcendental entire function of finite order satisfying 
	\begin{equation*}\label{LWC}
	T(r,A) \sim \log M(r,A)
	\end{equation*}
as $r\to\infty$ outside an exceptional set $G$ of finite logarithmic measure. If $E$ is a product of two linearly independent solutions of \eqref{ode}, then $\lambda(E)=\infty$.
\end{thm}

We prove the following generalization of Theorem~\ref{laine-wu}, which also improves the inequality \eqref{bl-inequality} when $1\leq \rho(A)< \frac{1-\beta}{2(1-\alpha)}$.

\begin{theorem}\label{thm}
Let $\alpha \in (0,1]$, and let $A$ be a transcendental entire function satisfying 
	\begin{equation}\label{regularity}
	T(r,A) \sim \alpha \log M(r,A) 
	\end{equation}
as $r\to\infty$ outside a set $G$ with $\underline{\operatorname{logdens}}(G)=\beta<1$. 
Suppose further that one of the following holds:
	\begin{center}
    \textnormal{(1)} $\rho(A)\notin \N, \qquad$ \textnormal{(2)} $\mu(A)< \rho(A),\qquad$ 
    \textnormal{(3)} $\rho(A) < \frac{1-\beta}{2(1-\alpha)}$.
    \end{center}
If $E$ is a product of two linearly independent solutions of \eqref{ode}, then 
	\[
	 \lambda(E) \ge \frac{1-\beta}{2(1-\alpha)}.
	\]
In particular, if $\alpha=1$, then $\lambda(E)=\infty$.
\end{theorem}

If $A$ is Mittag-Leffler's function of order $\rho\in (1/2, (2+\pi) /(2\pi))$, then $A$ satisfies \eqref{regularity} with $\alpha= \frac{1}{\pi \rho}$, see \cite[p.~19]{Hayman}. Such functions $A$ are examples of entire functions with the property $\rho(A) < \frac{1}{2(1-\alpha)}$. Examples of entire functions $A$ satisfying
(1) or (2) in  Theorem~\ref{thm} are standard. Each of the conditions (1)--(3) is necessary since the equation
	\[
	f'' + (e^z -1/16)f=0
	\]
has two linearly independent solutions $f_1$ and $f_2$ such that $\lambda(f_1f_2) = 0$, see \cite[p.~107]{L}.  Here the coefficient $A(z)=e^z -1/16$ has order $\rho(A)=1$ and satisfies \eqref{regularity} for $\alpha=1/\pi$. 

The lower bound of $\lambda(E)$ in Theorem~\ref{thm} does not depend on $\rho(A)$, and we can see that $\lambda(E)$ can be arbitrarily large provided only that $\alpha$ is close enough to $1$ independently
on the magnitude of $\rho(A)$.

\begin{corollary}\label{cor}
Let $A$ be a transcendental entire function, 
and let $\alpha = 1/ \beta^+(\infty,A)$. Suppose further that one of (1)--(3) with $\beta=0$ 
in Theorem~\ref{thm} holds. If $E$ is a product of two linearly independent solutions 
of \eqref{ode}, then 
	\[
	 \lambda(E) \ge \frac{1}{2(1-\alpha)}.
	\]
In particular, if $\alpha=1$, then $\lambda(E)=\infty$.
\end{corollary}


\section{Growth of solutions}\label{growth-sec}

It is well known that the coefficients $A$ and $B$ of \eqref{ode1} are polynomials if and only if
the solutions of \eqref{ode1} are of finite order. The possible orders in terms of the
degrees of $A$ and $B$ can be found in \cite{GSW}. Hence, if $A$ is transcendental and 
if $f_1$ and $f_2$ are linearly independent solutions of
\eqref{ode1}, then at least one of them is of infinite order. Finite order solutions are
also possible -- for example, $f(z)=e^{-z}$ solves \eqref{ode1} with $A(z)=e^z$ and
$B(z)=e^z-1$. This background led to asking the following research question
in \cite{G1}: \emph{What conditions on $A$ and $B$ will guarantee that every solution
$f\not\equiv 0$ of \eqref{ode1} has infinite order?} Examples of such conditions are
\begin{itemize}
\item[(i)] $\rho(A)<\rho(B)$, 
\item[(ii)] $A$ is a polynomial and $B$ is transcendental, 
\item[(iii)] $\rho(B)<\rho(A)\leq 1/2$,
\item[(iv)] $A$ is transcendental with $\rho(A)=0$ and $B$ is a polynomial,
\end{itemize}
see Theorems~2 and 6 in \cite{G1} and the main result in \cite{HMR}. 	

The seminal paper \cite{G1} has prompted a considerable amount interest in studying the growth of 
solutions of complex linear differential equations having well over one hundred citations
in the MathSciNet database in 2020. The following result by Laine and Wu is of particular
interest from the point of view of this paper.

\begin{thm}[\cite{LW2}]\label{laine-wu2}
Suppose that $A$ and $B\not\equiv 0$ are entire functions such that $\rho(B)<\rho(A)<\infty$ and
	$$
	T(r,A)\sim \log M(r,A)
	$$
as $r\to\infty$ outside a set $G$ of finite logarithmic measure. Then every non-trivial solution of \eqref{ode1} is of infinite order.
\end{thm}

Kwon and Kim \cite{KK} have showed that the conclusion of Theorem~\ref{laine-wu2} still holds if
the set $G$ satisfies $\overline{\operatorname{logdens}}(G)<(\rho(A)-\rho(B))/\rho(A)$.
If the condition $\rho(B)<\rho(A)<\infty$ is replaced with $\mu(B)<\mu(A)<\infty$, where $A$ satisfies
\eqref{regularity} as $r\to\infty$ outside a set $G$ satisfying $\overline{\operatorname{logdens}}(G)=0$,
then \cite[Theorem~1.5]{LHY} shows that
	$$
	\rho(f)\geq \frac{\mu(A)-\mu(B)}{21(\mu(A)+\mu(B))\sqrt{2\pi(1-\alpha)}}-1
	$$
for every solution $f\not\equiv 0$ of \eqref{ode1}. In particular, if $\alpha=1$, then $\rho(f)=\infty$.

For an entire function $g$, we define
	$$
	\xi(g) := \frac{1}{2\pi}\cdot \mes\left(\left\{\theta \in [0,2\pi) : \limsup_{r\to\infty} \frac{\log^+ |g(re^{i\theta})|}{\log r}< \infty \right\}\right). 
	$$
Clearly $0\le \xi(g)\le 1$. 
We have $\xi(g)=1$ if $g$ is a polynomial, and $\xi(g)=0$ if $g(z)=e^z+e^{-z}$. 
A transcendental entire function $g$ with $\xi(g)=1$ exists, see \cite[Lemma~4.1]{Hayman}. 
If $g$ is a Mittag-Leffler's function of order $\rho > 1/2$, then $\xi(g) = 1-\frac{1}{2\rho}$, see \cite[p.~19]{Hayman}.

The following result gives a new condition on the coefficients of \eqref{ode1} in terms of Petrenko's
deviation forcing the solutions to be of infinite order.

\begin{theorem}\label{g}
Let $A$ be an entire function such that $\xi(A)>0$, and let $B$ be a transcendental entire function satisfying $\beta^-(\infty,B)< \frac{1}{1 - \xi(A)}$. Then every non-trivial solution of \eqref{ode1} is of infinite order.
\end{theorem}

It follows from \eqref{Fabry} that an entire function $g$ with Fabry gaps satisfies $\beta^-(\infty,g)=1$. This gives raise to the following immediate consequence of Theorem~\ref{g}.

\begin{corollary}\label{coro}
Let $A$ and $B$ be entire functions. Suppose there exists a sector where $\log^+|A(z)| \lesssim \log|z|$, and suppose that $B$ is transcendental with Fabry gaps. Then every non-trivial solution of \eqref{ode1} is of infinite order.
\end{corollary}

Corollary~\ref{coro} improves \cite[Theorem~1.3]{Long} in the case when $A$ is a shortage
solution of $w''+P(z)w=0$ for a non-constant polynomial $P$ \cite{G0}. Indeed, for such $A$ there is
a sector in which $A$ tends to zero exponentially. This is particularly true if $A$ is the Airy
integral $\operatorname{Ai}(z)$ that solves the equation $w''-zw=0$. More generally, if $A\not\equiv 0$ is a contour integral
solution of 
	$$
	w^{(n)}+(-1)^{n+1}bw^{(k)}+(-1)^{n+1}zw=0,\quad n\geq 2,\ n>k>0,\ b\in\C,
	$$
then \cite[Theorem~3]{GHW} reveals that $\xi(A)\geq \frac{1}{2\pi}\cdot\frac{n\pi}{n+1}\geq \frac{1}{3}>0$. 

Corollary~\ref{coro} also improves \cite[Theorem~{1.7}]{Long}. Indeed, if $A$ is extremal for Yang's inequality, that is, if $p=q/2$, where $p$ denotes the number of finite deficient values and $q$ denotes the number of Borel's directions of order of $\ge \mu(A)$ of $A$, then \cite[Theorem~4]{W} asserts that there exists a sector where $A$ decays to a certain value $a\in\C$.

Using the $\cos\pi \rho$ -theorem, one can easily see that if $\xi(A)>0$ and $\mu(B)<1/2$, then every non-trivial solution of \eqref{ode1} is of infinite order.  The same conclusion holds if 	
	\begin{equation}\label{condition}
	1/2 \le \mu(B) < \frac{1}{\pi(1-\xi(A))}.
	\end{equation} 
This follows by \eqref{petrenko} and Theorem~\ref{g}.  Next we show that the condition \eqref{condition} can be weakened to $\mu(B)<\frac{1}{2(1-\xi(A))}$.

\begin{theorem}\label{h}
Let $A$ be an entire function such that $\xi(A)>0$, and let $B$ be a transcendental entire function satisfying $\mu(B)< \frac{1}{2(1-\xi(A))}$. Then every non-trivial solution of \eqref{ode1} is of infinite order.
\end{theorem}

To illustrate this theorem, let $A$ be Mittag-Leffler's function of order $\rho(A)>1/2$, and let $B$ be a transcendental entire function with $\mu(B) \neq \rho(A)$. Then $\xi(A) =1- \frac{1}{2\rho(A)}$ \cite[p.~19]{Hayman}, so that either $\mu(B)<\frac{1}{2(1-\xi(A))}$ or $\rho(B)\geq\mu(B)>\rho(A)$. It follows from Theorem~\ref{h} and \cite[Corollary~1]{G1} that every non-trivial solution of \eqref{ode1} is of infinite order.


\section{Proofs of Theorem~\ref{thm} and Corollary~\ref{cor}}\label{proof1-sec}

Let $g$ be an entire function, and let $D=\{z\in\C: |g(z)|>1\}$. For any $r>0$, let $\mathcal{A}_k(r)$ for  $k=1,2, \ldots, n(r)$ be the arcs of $|z|=r$ contained in $D$, and let $r\theta_k(r)$ be their lengths. Define $\theta(r)=\infty$ if the entire circle $|z|=r$ lies in $D$. Otherwise, define $\theta(r)=\max_k \theta_k(r)$. 

\begin{lemma}[\cite{A}]\label{arima} 
For any entire function $g$ and for any $0<\eta<1$, we have
	\begin{equation}\label{inequality-arima}
	\log\log M(r,g)> \pi \int_{r_0}^{\eta r} \frac{dt}{t \theta(t)} - c(\mu,r_0),
	\end{equation}
where $0<r_0<\eta r$ and $c(\mu, r_0)$ is a constant independent of $r$.
\end{lemma}

\begin{lemma}[{\cite{G}}]\label{log-derivative}
Let $g$ be a meromorphic function of finite order $\varrho$, and let $\varepsilon>0$ be a given constant. Then there exists a set $\F\subset [1,\infty)$ of finite logarithmic measure, such that for all $z$ satisfying $|z|\notin \F \cup [0,1]$ and for all integers $k>j\ge 0$, we have
	\[
	\left| \frac{g^{(k)}(z)}{g^{(j)}(z)} \right| < |z|^{(k-j)(\varrho-1+\varepsilon)}.
	\]
\end{lemma}

\begin{lemma}[{\cite{T}}] \label{toda}
Let $A$ and $E$ be entire functions satisfying \eqref{A} below. Suppose that $\lambda(E) < \rho(E)$. Then  
	$$
	\mu(E) = \rho(E) = \mu(A) = \rho(A),
	$$ 
and these numbers are equal to an integer or $\infty$.
\end{lemma}

We proceed to prove Theorem~\ref{thm} by modifying the reasoning in \cite{LW}.
Let $f_1$ and $f_2$ be two linearly independent solutions of \eqref{ode}, and set $E=f_1f_2$. 
From the results of Bank-Laine, Rossi and Shen, if $\rho(A):=\rho\le 1/2$, then $\lambda(E)=\infty$. Therefore, we may assume that $\rho>1/2$. We make use of the famous Bank-Laine formula \cite{BL2}
	\begin{equation}\label{A}
	-4A(z)= \frac{c^2}{E^2}+2\frac{E''}{E}-\left( \frac{E'}{E} \right)^2,
	\end{equation}
where $c$ is non-zero constant. Hence, 
	\begin{equation} \label{TrE}
	T(r,E)= N(r,1/E) + 2^{-1}T(r,A)+ S(r,E), \quad r\to\infty.
	\end{equation}
It follows from \eqref{TrE} that $\rho(E)$ and $\lambda(E)$ are both finite or both infinite. If $\rho(E)=\infty$, then there is nothing to prove, and for that reason we suppose that $\rho(E)=\varrho<\infty$. 

From Lemma~\ref{log-derivative}, choosing $\varepsilon=\frac{1}{2}$, we have
	\begin{equation}\label{pr1}
	2\left| \frac{E''(z)}{E(z)}\right| + \left| \frac{E'(z)}{E(z)} \right|^2 < 3|z|^{2\varrho-1},
	\end{equation}
for all $z$ satisfying $|z|\notin \F \cup [0,1]$, where $\F$ is a set of finite logarithmic measure. 
Set
	$$
	D_1:= \left\{z\in\C : |E(z)|>1 \right\}\quad\textnormal{and}\quad
	\F^* := \{z\in\C : |z|\in \F\},
	$$	
and let $r\theta_1(r)$ be the length of the longest arc of $|z|=r$ in $D_1$. 
Hence, from \eqref{A} and \eqref{pr1}, there is a constant $r_0>1$, such that for all $z\in D_1 \setminus \F^*$, $|z|>r_0$, we have 
	\begin{equation}\label{pr2}
	|A(z)|<|z|^{2\varrho}.
	\end{equation}
Define the sets 
	$$
	D_2:= \{z\in\C : |A(z)|>|z|^{2\varrho}\}
	\quad\textnormal{and}\quad
	H(r) := \{\theta \in [0,2\pi) : re^{i\theta} \in D_2 \}.
	$$
From \eqref{pr2}, it's clear that 
	\begin{equation}\label{D1D2}
	\Big( D_1\setminus (\F^* \cup \{|z|\le r_0\}) \Big) 
	\bigcap \Big( D_2\setminus (\F^* \cup \{|z|\le r_0\}) \Big)= \emptyset.
	\end{equation}
Then
	\begin{align*}
	2\pi T(r,A) &= \int_{H(r)} \log^+|A(re^{i\theta})| d\theta + \int_{\cm H(r)} \log^+|A(re^{i\theta})| d\theta \\
			&\le \mes(H(r)) \log M(r,A) + 2 \varrho (2\pi - \mes(H(r))) \log r,
	\end{align*}
which gives
	$$
	2\pi \le \mes(H(r)) \frac{\log M(r,A)}{T(r,A)} + 2 \varrho (2\pi - \mes(H(r))) \frac{\log r}{T(r,A)}.
	$$
Since $A$ is transcendental and satisfies \eqref{regularity} outside $G$, we obtain from the latter inequality that 
	\begin{equation}\label{mes}
	\liminf_{\substack{r\to\infty \\ r\notin G}} \mes(H(r)) \ge 2\pi \alpha.
	\end{equation}
Given $\veps>0$, from \eqref{D1D2} and \eqref{mes},  there exists $r_1>r_0$, such that for all $r \notin G \cup \F \cup [0,r_1]$, we have
	\begin{equation}\label{theta}
	\theta_1(r) \le \mes \left(\Big\{\theta\in [0,2\pi) : re^{i\theta} \in D_1 \Big\}\right) 
	\le 2\left(1-\alpha \right)\pi +\varepsilon.
	\end{equation}
	
Set $J_r := [r_1, r/2 ] \setminus \Big(G\cup\F \Big)$. Then it follows from Lemma~\ref{arima} and \eqref{theta}, that
	\begin{equation*}\label{arima_2}
	\log \log M(r,E) \ge \pi \int_{r_1}^{r/2} \frac{dt}{t\theta_1(t)} - c(1/2,r_1) \ge \frac{\pi}{2\left(1-\alpha \right)\pi+\veps} \int_{J_r} \frac{dt}{t}- c(1/2,r_1).
	\end{equation*}
Hence,
	\begin{eqnarray*}
	\rho(E) &\ge& \frac{\pi}{\left(1-\alpha \right) 2\pi + \veps} \Big(1- \lld(G\cup\F)\Big) \\
	&\ge& \frac{\pi}{\left(1-\alpha \right) 2\pi + \veps} \Big(1- \lld(G)-\uld(\F)\Big)
	=\frac{(1-\beta)\pi}{\left(1-\alpha \right) 2\pi + \veps}.
	\end{eqnarray*}
Letting $\veps \to 0^+$,   we get $\rho(E)=\infty$ if $\alpha=1$ and $\rho(E) \ge \frac{1-\beta}{2 \left(1-\alpha \right)}$ if $\alpha<1$.
Since $A$ satisfies one of (1)--(3), it follows from Lemma~\ref{toda} that $\rho(E)=\lambda(E)$,
which in turn implies the conclusion of Theorem~\ref{thm}.

\medskip
It remains to prove Corollary~\ref{cor}. We argue similarly as above up to \eqref{mes}, which 
now reads without the exceptional set $G$.
If $\alpha>0$, we get \eqref{theta}, which yields the lower bound $\rho(E)\geq\frac{\pi}{\left(1-\alpha \right) 2\pi + \veps}$, where we may let $\veps\to 0^+$. If $\alpha=0$, that is, if $\beta^+(\infty,A)=\infty$,
we use $\theta_1(r)\leq 2\pi$, and obtain $\lambda(E)\geq\frac12$. Finally we apply Lemma~\ref{toda}.


\section{Proof of Theorems~\ref{g} and \ref{h}}\label{proof2-sec}

We recall the following lemmas.

\begin{lemma}[{\cite{G}}]\label{log-derivative-2}
Let $g$ be a meromorphic function of finite order $\varrho$, and let $\varepsilon>0$ be a given constant. Then there exists a set $E\subset[0,2\pi)$ that has linear measure zero, such that if $\psi_0 \in [0,2\pi)\setminus E$, then there exists a constant $R_0=R_0(\psi_0)>1$ such that for all $z$ satisfying $\arg z=\psi_0$ and $|z|\ge R_0$, and for all integers $k>j\ge 0$, we have
	\[
	\left| \frac{g^{(k)}(z)}{g^{(j)}(z)} \right| < |z|^{(k-j)(\varrho-1+\varepsilon)}.
	\]
\end{lemma}

\begin{lemma}[\cite{WLHQ}]\label{Ph-Li}
Let $g$ be an entire function of lower order $\mu(g)\in [1/2,\infty)$. Then there exists a sector $S(\alpha,\beta) = \{z: \alpha<\arg z <\beta\}$ with $\beta-\alpha> \frac{\pi}{\mu(g)}$ and $0\le \alpha<\beta\le 2\pi$, such that 
	$$
	\limsup _{r \rightarrow \infty} \frac{\log \log \left|g\left(r e^{i \theta}\right)\right|}{\log r} \geq \mu(g)
	$$
holds for all rays $\arg z = \theta \in (\alpha,\beta)$.
\end{lemma}

We proceed to prove Theorem~\ref{g}.
Suppose on the contrary to the assertion that there exists a non-trivial solution $f$ of \eqref{ode1} with $\rho(f)=\r < \infty$. Then, from Lemma~\ref{log-derivative-2}, we have for any $\theta\in [0,2\pi) \setminus E$, where $E \in [0,2\pi)$ has linear measure zero, that there is a constant $R(\theta)>1$ such that for any $r> R(\theta)$, 
	\begin{equation}\label{G1-1}
	\left| \frac{f^{(j)}(re^{i\theta})}{f(re^{i\theta})} \right| \le r^{2\r}, \quad j=1,2.
	\end{equation}
Recall that $1\le\beta^-(\infty,B)<1/(1-\xi(A))$, where $\xi(A)>0$. Given constants 
	$$
	0< \veps < \frac{1}{\beta^-(\infty,B)}- (1-\xi(A))
	\quad\textnormal{and}\quad
	\frac{2}{2+ \veps}<d<1,
	$$ 
define
	$$
	I_d(r) := \{\theta \in [0,2\pi) : \log |B(re^{i\theta})| \ge  (1-d) \log M(r,b)\}.
	$$
Then, 
	\begin{equation}\label{mes2}
	\begin{split}
	2\pi T(r,B) &= \int_{I_d(r)} \log^+ |B(re^{i\theta})| d\theta + \int_{\cm I_d(r)} \log^+ |B(re^{i\theta})| d\theta\\
	& \le \mes(I_d(r)) \log M(r,B)  + \Big(2\pi - \mes(I_d(r))\Big) (1-d) \log M(r,B).
	\end{split}
	\end{equation}
Dividing both sides of \eqref{mes2} by $\log M(r,B)$ and using the definition \eqref{deviation}, we deduce
	\begin{equation}\label{G1-3}
	\limsup_{r\to \infty} \mes(I_d(r)) \ge \frac{2\pi}{d \beta^-(\infty,B)} - \frac{2\pi (1-d)}{d}.
	\end{equation}
For the choice of $\veps$ and $d$, we deduce from \eqref{G1-3}, that there exist an infinite sequence $\{r_\nu\}$ and $R^*>0$, such that for all $\nu \in \N$ for which $r_\nu>R^*$, we have
	\begin{equation}\label{G1-4}
	\mes(I_d(r_\nu)) \ge \frac{2\pi}{d \beta^-(\infty,B)} - \frac{2\pi (1-d)}{d} -\pi \veps 
	>2\pi(1 - \xi(A)).
	\end{equation}
Thus there exists an interval $(\theta_1, \theta_2)$ such that
	$$
	(\theta_1, \theta_2) \subset I_d (r_\nu) \cap \left\{\theta \in [0,2\pi) : \limsup_{r\to\infty} \frac{\log^+ |A(re^{i\theta})|}{\log r}< \infty \right\}.
	$$ 
Therefore, for any $\theta \in (\theta_1, \theta_2) \setminus E$, we obtain by using \eqref{ode1} and \eqref{G1-1},
	\begin{align*}
	\log^+ M(r_\nu, B) &\lesssim \log^+ \left|B(r_\nu e^{i\theta})\right| \\
					& \lesssim \log^+ \left| \frac{f''(r_\nu e^{i\theta})}{f(r_\nu e^{i\theta})}\right| + \log^+ \left| \frac{f'(r_\nu e^{i\theta})}{f(r_\nu e^{i\theta})}\right| + \log^+ \left|A(r_\nu e^{i\theta})\right| + 1\\
					& \lesssim \log r_\nu, \quad  \nu \to \infty.
	\end{align*}
This implies that $B$ is a polynomial, which contradicts the assumption that $B$ is transcendental. Thus, every non-trivial solution of \eqref{ode1} is of infinite order.	

\bigskip

We proceed to prove Theorem~\ref{h}. 
If $\mu(B)< 1/2$, then by using the $\cos \pi \rho$ -theorem, we get the conclusion of the theorem. Hence we assume that $\frac{1}{2} \le \mu(B)< \frac{1}{2(1-\xi(A))}$.
Suppose on the contrary to the assertion that there exists a non-trivial solution $f$ of \eqref{ode1} with $\rho(f)=\r < \infty$. Then, from Lemma~\ref{log-derivative-2}, we have for any $\theta\in [0,2\pi) \setminus E$, where $E \in [0,2\pi)$ has linear measure zero, that there is a constant $R(\theta)>1$ such that for any $r> R(\theta)$,  \eqref{G1-1} holds. From Lemma~\ref{Ph-Li}, there is a sector $S(\alpha,\beta)= \{z: \alpha< \arg z<\beta\}$ with 
	$$
	\beta-\alpha \ge \frac{\pi}{\mu(B)} > 2\pi(1 - \xi(A))
	$$
such that 
	\begin{equation*}
	\limsup _{r \rightarrow \infty} \frac{\log \log \left|B\left(r e^{i \theta}\right)\right|}{\log r} \geq \mu(B)
	\end{equation*}
holds for all rays $\arg z = \theta \in (\alpha,\beta)$. Thus there exists an interval $(\theta_1, \theta_2)$ such that
	$$
	(\theta_1, \theta_2) \subset ( \alpha,\beta) \cap \left\{\theta \in [0,2\pi) : \limsup_{r\to\infty} \frac{\log^+ |A(re^{i\theta})|}{\log r}< \infty \right\}.
	$$ 
Therefore, there exists a sequence $z_n=r_n e^{i\theta}$ with $r_n \to \infty$ as $n\to\infty$ and $\theta\in (\theta_1,\theta_2)\setminus E$ such that 
	$$
	\begin{aligned}
\exp \left(r_{n}^{\mu(B)-\varepsilon}\right) & \leq |B(r_{n} e^{i \theta})| \\
& \leq\left|\frac{f^{\prime \prime}\left(r_{n} e^{i \theta}\right)}{f\left(r_{n} e^{i \theta}\right)}\right|+|A(r_{n} e^{i \theta})|\left|\frac{f^{\prime}\left(r_{n} e^{i \theta}\right)}{f\left(r_{n} e^{i \theta}\right)}\right| \\
& \leq r_{n}^{2 \r}(1+o(1)),
\end{aligned}
	$$
where $n$ is large enough and $\varepsilon>0$ is small. But this is a contradiction, and so $\rho(f)=\infty$ for every non-trivial solution $f$ of \eqref{ode1}.


\textsc{Address:} University of Eastern Finland, Department of Physics and Mathematics,
P.~O.~Box 111, 80100 Joensuu, Finland

\textsc{E-mail:} \texttt{janne.heittokangas@uef.fi, amine.zemirni@uef.fi}

\begin{thebibliography}{99}
\bibitem{A} K.~Arima,
\emph{On maximum modulus of integral functions},
J.~Math.~Soc.~Japan \textbf{4} (1952), 62--66. 

\bibitem{BL} S.~B.~Bank and I.~Laine, 
\emph{On the oscillation theory of $f''+Af = 0$ where $A$ is entire},
Bull.~Amer.~Math.~Soc.~(N.S.) \textbf{6} (1982), no.~1, 95--98.

\bibitem{BL2} S.~B.~Bank and I.~Laine, 
\emph{On the oscillation theory of $f''+Af = 0$ where $A$ is entire},
Trans.~Amer.~Math.~Soc.~\textbf{273} (1982), no.~1, 351--363.


\bibitem{BE2} W.~Bergweiler and A.~Eremenko, 
\emph{Quasiconformal surgery and linear differential equations},
J.~Anal.~Math.~\textbf{137} (2019), no.~2, 751--812.

\bibitem{BE} W.~Bergweiler and A.~Eremenko, 
\emph{On the Bank--Laine conjecture},
J.~Eur.~Math.~Soc.~(JEMS) \textbf{19} (2017), no.~6, 1899--1909.

\bibitem{C} J.~Clunie, 
\emph{On integral functions having prescribed asymptotic growth},
Canadian J.~Math.~\textbf{17} (1965) 396--404.

\bibitem{F} W.~H.~J.~Fuchs,
\emph{Proof of a conjecture of G. P\'olya concerning gap series},
Illinois J.~Math.~\textbf{7} (1963), 661--667. 



\bibitem{G}  G.~G.~Gundersen,
\emph{Estimates for the logarithmic derivative of a meromorphic function, plus similar estimates}, 
J.~London Math.~Soc.~(2) \textbf{37} (1988), 88--104.

\bibitem{G1} G.~G.~Gundersen,
\emph{Finite order solutions of second order linear differential equations},
Trans.~Amer.~Math.~\textbf{305} (1988), no.~1, 415--429.

\bibitem{G0} G.~G.~Gundersen,
\emph{On the real zeros of solutions of $f''+A(z)f=0$ where $A(z)$ is entire},
Ann.~Acad.~Sci.~Fenn.~Ser.~A I Math.~\textbf{11} (1986), no.~2, 275--294.

\bibitem{GHW} G.~G.~Gundersen, J.~Heittokangas and Z.-T.~Wen,
\emph{Families of solutions of differential equations that are defined by contour integrals},
\url{https://arxiv.org/abs/1911.09479}, 30 p.

\bibitem{G2} G.~G.~Gundersen and E.~M.~Steinbart,
\emph{A generalization of the Airy integral for $f''-z^nf=0$},
Trans.~Amer.~Math.~\textbf{337} (1993), no.~2, 737--755.

\bibitem{GSW}
G.~G.~Gundersen, E.~Steinbart and S.~Wang, 
\emph{The possible orders of solutions of linear differential equations with polynomial coefficients}, 
Trans.~Amer.~Math.~Soc.~\textbf{350} (1998), no.~3, 1225--1247.

\bibitem{Hayman} W.~K.~Hayman,
\emph{Meromorphic Functions}.
Oxford Mathematical Monographs, Clarendon Press, Oxford, 1964.

\bibitem{HR} W.~K.~Hayman and J.~F.~Rossi,
\emph{Characteristic, maximum modulus and value distribution},
Trans.~Amer.~Math.~\textbf{284} (1984), no.~2, 651--664.

\bibitem{HMR} S.~Hellerstein, J.~Miles and J.~Rossi,
\emph{On the growth of solutions of $f''+gf'+hf=0$},
Trans.~Amer.~Math.~\textbf{324} (1991), no.~2, 693--706.

\bibitem{KK} K.~Kwon and J.~Kim,
\emph{Maximum modulus, characteristic, deficiency and growth 
of solutions of second order linear differential equation},
Kodai Math.~J.~\textbf{24} (2001), no.~3, 344--351.

\bibitem{L} I.~Laine,
\emph{Nevanlinna Theory and Complex Differential Equations}. 
De Gruyter Studies in Mathematics, 15. Walter de Gruyter \& Co., Berlin, 1993.

\bibitem{LW2} I.~Laine and P.~C.~Wu, 
\emph{Growth of solutions of second order linear differential equations}, 
Proc.~Amer.~Math.~Soc.~\textbf{128} (2000), no.~9, 2693--2703.

\bibitem{LW} I.~Laine and P.~C.~Wu,
\emph{On the oscillation of certain second order linear differential equations}, 
Rev.~Roumaine Math.~Pures Appl.~\textbf{44} (1999), no.~4, 609--615.

\bibitem{Long} J.~Long,
\emph{Growth of solutions of second order complex linear differential equations with entire coefficients},
Filomat \textbf{32} (2018), no.~1, 275--284.

\bibitem{LHY} J.~Long, J.~Heittokangas and Z.~Ye,
\emph{On the relationship between the lower order of coefficients and the growth 
of solutions of differential equations},
J.~Math.~Anal.~Appl.~\textbf{444} (2016), no.~1, 153--166.


\bibitem{P} V.~P.~Petrenko, 
\emph{Growth of meromorphic functions of finite lower order},  
Izv.~Akad.~Nauk SSSR Ser.~Mat.~\textbf{33} (1969), 414--454 (Russian). 
English transl.: Math.~USSR, Izv.~\textbf{3} (1969), 391--432.

\bibitem{R} J.~Rossi,
\emph{Second order differential equations with transcendental coefficients},
Proc.~Amer.~Math.~Soc.~\textbf{97} (1986), no.~1, 61--66.


\bibitem{Shen} L.~C.~Shen,
\emph{Solution to a problem of S.~Bank regarding exponent of convergence of zeros of the solutions of differential equation $f''+Af=0$}, 
Kexue Tongbao (English ed.) \textbf{30} (1985), no.~12, 1579--1585.


\bibitem{T} N.~Toda, 
\emph{A theorem on the growth of entire functions on asymptotic paths and its application to the oscillation theory of $w'' + Aw = 0$},
Kodai Math.~J.~\textbf{16} (1993), no.~3, 428--440.


\bibitem{W} S.~J.~Wu, 
\emph{Some results on entire functions of finite lower order},
Acta Math.~Sinica (N.S.) \textbf{10} (1994), no.~2, 168--178.



\bibitem{WLHQ}
X.~Wu, J.~Long, J.~Heittokangas and K.~E.~Qiu,
\emph{On Second-order complex linear differential equations with special functions or extremal functions as coefficients}, 
Electron. J. Differential Equations {\bf 2015}, No. 143, 15 pp.

\end{thebibliography}
\end{document}